\documentclass[12pt]{article}
\usepackage{amssymb}
\usepackage{amsmath}

\begin{document}

\begin{large}

\begin{center}
\begin{LARGE}
\vspace*{\fill}
\vspace*{\fill}
\vspace*{\fill}

{\bf THE BLOCK RELATION IN COMPUTABLE LINEAR ORDERS}
\vspace*{\fill}
\vspace*{\fill}

\end{LARGE}

{\bf Michael Moses} \\
{\bf moses@gwu.edu}
\vspace*{\fill}

{\bf Department of Mathematics} \\
{\bf The George Washington University}
\end{center}
\vspace*{\fill}
\vspace*{\fill}

\noindent
{\it 2000 Mathematics Subject Classification}: 03D45 (primary), 03C57 (secondary) \\
{\it Keywords and phrases}: computable linear order, block relation, self-embedding 

\newpage
\begin{center}
{\bf Abstract}
\end{center}

\bigskip

\noindent

A block in a linear order is an equivalence class when factored by the block relation $B(x,y)$, satisfied by elements that are finitely far apart. We show that every computable linear order with dense condensation-type (i.e.\ a dense collection of blocks) but no infinite, strongly $\eta$-like interval (i.e.\ with all blocks of size less than some fixed, finite $k$) has a computable copy with the non-block relation $\neg B(x,y)$ computably enumerable. This implies that every computable linear order has a computable copy with a computable non-trivial self-embedding, and that the long-standing conjecture characterizing those computable linear orders every computable copy of which has a computable non-trivial self-embedding (as precisely those that contain an infinite, strongly $\eta$-like interval) holds for all linear orders with dense condensation-type. 

\newpage
\begin{center}
{\bf Introduction}
\end{center}

\bigskip

I have always found attractive those several results of mathematical logic that establish an equivalence between the (syntactic) definability of a mathematical property in a particular mathematical language and the (semantic) `algebraic' characteristics of that property. Representative of these results, and one which has provided a focus for much of my research, is the famous result of Ash and Nerode [1981] in which (under certain `extra decidability' conditions) a relation is shown to be {\it intrinsically} computable in a computable mathematical structure (meaning that it is computable in every computable copy of that structure) if and only if it and its complement are definable in the structure by computable $\Sigma_{1}$ $\mathcal{L}_{\omega_{1}\omega}$ formula). Characterizing the other side of the coin, intrinsically {\it non}-computable relations (for which {\it no} computable copy of the structure exists in which the relation is computable) will prove more difficult; the defining formulae required will be complex. Recently developed tools such as the `separator' construction of Jockusch and Soare [1991] and other techniques used to analyze the degree spectra of relations and structures may suggest approaches that were not evident before. The main result of this paper was arrived at as part of a (continuing) project along those lines, seeking to characterize those isomorphism types in a particular theory (linear order) in which a particular relation (the block relation) is intrinsically non-computable, but is presented here outside of that context because of some interesting ramifications with regard to computable non-trivial self-embeddings of computable linear orders. 

\medskip

A computable linear order $L$ is one with universe $\mathbb{N}$ on which the order relation $<_{L}$ is computable. The block relation $B_{L}(x,y)$ is that satisfied by elements $x,y$ for which the $L$-interval $[x,y]$ is finite. In one of the earliest results about these structures, Tennenbaum showed (see Rosenstein [1982]) that there is a computable linear order of type $\omega + \omega^{\star}$ in which the block relation is not computable. In [1970] Feiner established the existence of a computable linear order no computable copy of which has the block relation computable; and in [1986] Moses characterized intrinsically computable relations on linear orders (showing that the Ash-Nerode characterization holds for linear orders without the need for `extra decidability' conditions).

\medskip

The main result of this paper establishes that every computable linear order with dense condensation-type but no infinite, strongly $\eta$-like interval has a computable copy with computably enumerable non-block relation. The {\it condensation-type} of a linear order is the order-type of its `finite condensation', i.e.\ the linear order obtained by factoring out by the block relation, clumping each block into a single point. A {\it strongly $\eta$-like} interval is one (with condensation-type $\eta$) all of whose blocks are smaller than some fixed, finite $k$. 

\medskip

In [1940] Dushnik and Miller observed that every computable linear order has a non-trivial self-embedding. (That paper is famous for a non-computable `finite injury' construction of a dense suborder of $\mathbb{R}$ with no non-trivial self-embedding.) The characterization of linear orders every computable copy of which has a computable non-trivial self-embedding remains open, despite repeated salvos, and despite the straightforward and long-standing conjecture (Kierstead [1987]): Every computable copy of a computable linear order has a computable non-trivial self-embedding if and only if the linear order contains an infinite, strongly $\eta$-like interval. A corollary to the main result of this paper establishes that this conjecture holds for linear orders whose condensation-type is dense. (This extends slightly a recent result of Downey, Kastermans, and Lempp [2009] in which the characterization is shown to hold for all $\eta$-like linear orders, i.e.\ those with dense condensation-type {\it and} with no infinite block.) A surprising, at least to me, second corollary establishes that {\it every} computable linear order has a computable copy with a computable non-trivial self-embedding.

\newpage
\begin{center}
{\bf The Main Result}
\end{center}

\bigskip

\noindent
{\bf THEOREM 1:}

\medskip

{\it Every computable linear order of condensation-type $\eta$ with no infinite, strongly $\eta$-like interval has a computable copy with computably enumerable non-block relation.} 

\medskip

\noindent
{\bf PROOF:}

\medskip

Let $L$ be a computable linear order of condensation-type $\eta$ with no infinite, strongly $\eta$-like interval; we shall construct a computable linear order $M$ with computably enumerable non-block relation and an isomorphism $f: L \to M$. 

\bigskip

\noindent
{\bf {\it Preliminary Definitions}}: 

\medskip

Consider the universe $\mathbb{N}$ of $L$ as being enumerated in stages, one element at a time, with $L^{s}$ denoting the finite linear order enumerated by stage $s$, on $\{ 1, 2, \cdots, s\}$, ordered according to $L$'s computable order relation $<_{L}$. We will define a computable linear order $M$ by enumerating its universe $\mathbb{N}$ in stages (usually several more than just one element at each stage), with order relation $<_{M}$ inherited from $L$, which order will never subsequently be changed, thus making $M$ a computable linear order. We will also enumerate at each stage some pairs of $M$-elements into $M$'s non-block relation $\neg B_{M}(x,y)$, which enumeration will never subsequently be withdrawn (thus making it computably enumerable). We will also define at each stage $s$ a partial function $f^{s}: L^{s} \to M^{s}$ which sequence of functions will define the isomorphism $f: L \to M$ as the limit along a subsequence of stages; i.e.\ we will establish the existence of an infinite subsequence of stages $s_{1}, s_{2}, s_{3}, \cdots$ such that, for each $a$ in $L$ and $b$ in $M$, $f^{s_{i}}$ will remain unchanged on $a$ and $b$ for all but finitely many $i$, which limiting values are defined to be $f(a)$ and $f^{-1}(b)$. 

\medskip

Our strategy to make $M$ isomorphic to $L$ while keeping $\neg B_{M}(x,y)$ computably enumerable in $M$ is straightforward: At each stage of our construction we will seek to define an image in $M$ for each $L$-block and enumerate pairs of $M$-elements into $\neg B_{M}(x,y)$ if they lie in images of different $L$-blocks. We will access those $L$-blocks via their least-block-elements, i.e.\ the $\mathbb{N}$-least element of each block in $L$, defined by the $\Pi_{2}$ formula $\forall y <_{\mathbb{N}} x \; \exists z \: ( (y <_{L} z <_{L} x) \vee (x <_{L} z <_{L} y))$. We will use a computable binary relation $R(x,y)$ that underlies this $\Pi_{2}$ formula (in the sense that $n$ is a least-block-element if and only if $R(n,y)$ holds for infinitely many $y$); an $L$-element $n$ appears to be a least-block-element at stage $s$ (we will say that $n$ is {\it on} at that stage) if $R(n,s)$ holds. The true least-block-elements in $L$, and only the true least-block-elements, will be {\it on} at infinitely many stages.  As observed in Jockusch [1966] we can select the $R(x,y)$ so that, for each $n$, there will be infinitely many stages at which the elements from among $1, 2, 3, \cdots, n$ that are {\it on} are precisely the true least-block-elements from among those first $n$ $L$-elements. Our construction will at each stage seek, for the $L$-elements $n$ that are {\it on}, an image in $M$ for $n$'s block in $L$, of the right size and in the right location. The image of an $L$-block whose least-block-element is not {\it on} will need to be incorporated into another $M$-cluster that is enumerated subsequently (so as not to compromise $\neg B_{M}(x,y)$), and incorporated whole, ready to come back to life when (if) that least-block-element next comes back {\it on}. The fact that every infinite interval in $L$ contains arbitrarily large blocks is just what we need in order to pull this off: a large enough true block will eventually appear where we need it in order to incorporate the image of an $L$-block whose least-block-element is no longer {\it on}. Some delicacy, nothing too onerous, is required to identify the true $L$-block whose image in $M$ will incorporate this fallow $M$-cluster. 

\medskip

The image in $M^{s}$ of the $L$-block around the least-block-element $n$ will carry an identifying label. This labeled cluster of elements will be contiguous in $M^{s}$, nor will an element ever be introduced internally at a later stage (so long as the cluster carries the label). Since we can guarantee only that a true least-block-element $n$ will appear {\it on} infinitely often (not that it will always appear {\it on} after some stage), we will need to preserve the image of its block through those stages when it is not {\it on}. During those stages we will seek to incorporate this $M$-cluster labeled for $n$ into images of other $L$-blocks, causing this cluster to take on additional labels, for some other seeming least-block-elements, which labels will be shaken off when $n$ next comes back {\it on}. Elements that are jettisoned in this way will never subsequently return to the $M$-cluster, thus maintaining the integrity of $M$'s non-block relation. So, for each true least-block-element $n$ of $L$, the image in $M$ of the $L$-block  around $n$ will (after a few false starts) develop as a contiguous cluster of labeled elements in $M$, which may contain within it another labeled cluster, denoting the fallow image in $M$ of the block  around an $L$-element that was once {\it on} but is no longer (which may contain another fallow block within it, etc.), and will, during those stages when $n$ is not {\it on}, become temporarily incorporated into other larger labeled clusters, which labels will subsequently be shaken off. This image for $n$'s block will grow and shrink according to the current size of $n$'s block in $L$, our approximation at each stage (which we define a couple of paragraphs below) to the final $L$-block around $n$, which is, unfortunately, no better than $\Sigma_{2}$-definable. We will, moreover, need to manage not just the size of the image: if the $L$-block around $n$ turns out to be infinite (i.e.\ one of $\omega$, $\omega^{\star}$ or $\omega^{\star} + \omega$), we must ensure that its image will be a block of the same type. 

\medskip

We will use a slightly souped-up version of the usual $\Pi_{2}$ guessing tree to drive our construction: The standard part imagined, as usual, as branching downward, with branching at level $n$ denoting the guess as to whether or not $n$ is indeed a true least-block-element in $L$, with the left branch denoting the positive guess and the right the negative. The (finite) path that is {\it on} at stage $s$ is the one from the root to the $s$th level that branches left at each level $n$ if $n$ is {\it on} at stage $s$ and branches right otherwise. The (infinite) {\it true} path is the one that branches left at each level $n$ precisely if $n$ is a true least-block-element in $L$. We have the usual features: for each level $n$, the {\it true} path is the leftmost one that coincides infinitely often with the {\it on} path at the first $n$ levels. We will drive a separate construction along each path of the tree, acting along the {\it on} path at each stage and using as labels on clusters in $M$ the node $\sigma$ at which those labels were created, which labeled $M$-cluster will denote the image of the $L$-block around the element that $\sigma$ references. As it is, however, this tree is too simple for our purposes: it will suffice to produce an image in $M$, of the right size and in the right location, for each $L$-block around a true least-block-element but will leave scattered throughout $M$ the fallow images of $L$-blocks centered around elements that once seemed to be least-block-elements but are no longer. In order to incorporate into these constructions the search for an appropriate preimage for each fallow $M$-cluster, i.e.\ an $L$-block centered around a true least-block-element into whose image we can incorporate the fallow $M$-cluster, we will expand our $\Pi_{2}$ guessing tree by adding immediately below each level $n$ two additional sublevels (described three paragraphs below), which levels will be dynamic, in that the $L$-elements referenced by those nodes will (may) change as the construction progresses.

\medskip

The block at stage $s$ around an $L$-element $n$ (a potential least-block-element) is defined to be the contiguous cluster of $L^{s}$-elements around $n$ extending as far as possible on both sides without including any element $\mathbb{N}$-less than $n$, nor any element that was enumerated into $L$ at or after the last stage (previous to $s$) when $n$ was {\it on} (or stage 1 if $n$ was never {\it on}), nor any element (other than $n$) that came {\it on} at or after the last stage when $n$ was {\it on}. Notice that, if $m$ is not a true least-block-element (which will cause it to be {\it on} at only finitely many stages), the $L$-block around $m$ will eventually stop changing, and will then be a sub-block of the (true) block around $m$. And, if $n$ is a true least-block-element (which will cause it to be {\it on} infinitely often) and $m$ is truly in $n$'s block, then, after finitely many stages, $m$ {\it and all of $m$'s block} will always be in $n$'s block (there will come a stage after which no new elements are enumerated between $m$ and $n$ and after which no elements between $m$ and $n$ are ever {\it on}). Notice also that, if $m$ and $n$ are truly in separate blocks in $L$, there will be infinitely many stages at which they will not be members of each other's blocks. Unfortunately, it could well happen that they also will be members of each other's blocks infinitely often; we must make sure that this does not inhibit the growth of the separate images in $M$ of their separate blocks in $L$. Finally note that our definition allows adjacent blocks in $L^{s}$ to overlap. 

\medskip

As mentioned previously, our main concern will be to incorporate every fallow $M$-cluster into the image of a true $L$-block. The fallow $M$-cluster around an $M$-element $n$ is the cluster defined by all the labels that $n$ carries of nodes on paths to the left of the node $\sigma$ at which we're acting, i.e.\ the cluster of elements around $n$ that also carry (any one of) those labels. Any preimage that we find for $n$ while acting at $\sigma$ must incorporate wholly, into a single $L$-block, all of this fallow $M$-cluster around $n$ (in anticipation of it coming back {\it on} at some later stage). If $n$ carries no such labels (of nodes left of $\sigma$), then there is no fallow $M$-cluster for us to worry about; however, rather than set up separate machinery to handle this simpler case, we will consider the singleton $n$, by itself, to be the fallow $M$-cluster around $n$. The search for a preimage for the fallow $M$-cluster around $n$ will begin rather simply: we will identify the size of the block that we need (the current size of the fallow $M$-cluster) and the $L$-interval (between the $L$-blocks around higher priority elements) in which we need it. We will maintain a lexicographically ordered list of contiguous sequences of $L$-elements of the required size within this interval, adding new contiguous sequences to the right end of the list as they are enumerated into $L$ and removing from the list those sequences that are no longer contiguous (within which a new $L$-element has been introduced). The leftmost on this list will be the selected preimage at that stage. Notice that, since every infinite interval in $L$ contains arbitrarily large blocks, the selected preimage will, after a finite number of stages, settle on a truly contiguous sequence of $L$-elements suitable for use as the preimage for the fallow $M$-cluster around $n$.

\medskip

Since our construction is centered around least-block-elements, we will need to identify the least-block-element in $L$ within whose block this selected preimage lies. We will maintain an ordered list of possible least-block-elements consisting of all $i$ that are $\leq_{\mathbb{N}}$ the $\mathbb{N}$-least element of the selected preimage (no other $i$ could possibly be the least-block-element whose block contains that preimage). We will order the $i$ on the list according to the last time that a new element was introduced between $i$ and the selected preimage, moving $i$ to the right end of the list each time a new element is introduced between $i$ and the selected preimage (ordering them as in $\mathbb{N}$ if there is more than one such $i$ being moved). Notice that, for a particular selected preimage, this procedure will eventually identify, as the leftmost element in this list that appears {\it on} infinitely often, the true least-block-element $i$ in $L$ whose block truly contains that selected preimage. (This $i$ will eventually cease being moved to the right end of the list and any elements to its left on the list will be in $i$'s $L$-block, i.e.\ the same block as the selected preimage, but will not be true least-block-elements and, consequently, will eventually cease appearing {\it on}. We manage this procedure by splitting our $\Pi_{2}$ guessing tree according to this list of all the $i$ that could possibly be the least-block-element whose $L$-block contains the selected preimage, which splitting will be dynamic in that the $i$ that these nodes will reference will (may) change from stage to stage. As noted, the $i$'s referenced by the leftmost of these nodes will eventually cease changing, and one of them, the true least-block-element whose block contains the selected preimage, will appear {\it on} infinitely often.  

\medskip

An obvious problem will occur if the true least-block-element $i$ so identified (whose block contains the selected preimage) turns out to be an $n$ of higher priority. The $L$-block around that $n$ will already have an image in $M$ and $M$'s non-block relation $\neg B_{M}(x,y)$ will not allow us to incorporate the fallow $M$-cluster into that existing image. So we need to make sure that the selected preimage for the fallow $M$-cluster is not contained within the blocks around the higher priority $n$. We do this by identifying a pair $p,q$ of true least-block-elements within the appropriate interval in $L$ (strictly between the higher priority $n$) and conducting our search for a preimage within the $L$-interval $[p,q]$. In order to settle on one such pair $p,q$ we will maintain a lexicographically ordered list of pairs $p,q$ of $L$-elements that lie within the appropriate interval between higher priority $n$, working at each stage within the interval defined by the leftmost pair on this list both of whose elements are {\it on}. Notice that, for a particular set of higher priority $n$, this procedure will identify, as the leftmost pair on this list that appears {\it on} infinitely often, a pair $p,q$ of true least-block-elements (within which infinite interval $[p,q]$ we are guaranteed to find a preimage that is not part of the block around any higher priority $n$). All the pairs to the left of $p,q$ will eventually cease appearing {\it on} and the preimage searches conducted between those pairs will be abandoned; all preimage searches conducted between pairs to the right of $p,q$ will be rescinded the next time $p,q$ comes {\it on}. We manage this procedure with another branching of our $\Pi_{2}$ guessing tree, with the nodes referencing in lexicographic order, from left to right, the $p,q$ in the abovementioned list, and conduct a search as described in the previous paragraph within each interval $[p,q]$ for a preimage for the fallow $M$-cluster and for the true least-block-element $i$ whose block contains the selected preimage. The pairs $p,q$ that these nodes reference will not change but new nodes will be added to the right end as new elements are enumerated into $L$. 

\medskip

This then is our $\Pi_{2}$ guessing tree, with three sublevels for each $n$: Branching at sublevel-one, if necessary, to manage the search for an image in $M$ for the $L$-block around the least-block-element $n$: Each node at the previous level will have a binary split at sublevel-one for $n$ denoting the guess as to whether $n$ is a true least-block-element in $L$ (left branch) or not (right). Branching at sublevels two and three, if necessary, to manage the search for a preimage in $L$ for the fallow $M$-cluster around $n$: Each of the nodes from sublevel-one will split at sublevel-two into finitely many branches denoting the lexicographically ordered list of all pairs $p,q$ that lie in $L^{s}$ between the higher priority elements between whose images the $M$-element $n$ lies, and within which interval $[p,q]$ we will seek to identify a preimage for the fallow $M$-cluster around $n$. Each of these sublevel-two nodes $p,q$ will split at sublevel-three into finitely many branches denoting the ordered list of possibilities for the least-block-element $i$ within $[p,q]$ whose block contains the selected preimage. Note that these second and third sublevels are dynamic in that the number of branches at sublevel-two will increase as more pairs $p,q$ become enumerated between the higher priority elements (although the $p,q$ that each node references will not change), and the $L$-elements $i$ referenced by the nodes at sublevel-three will change as elements are introduced into $L$ between $i$ and the selected preimage (although the number of branches will not change). Branching at these three sublevels for $n$ occurs only if necessary: If $n$ in $L$ is already referenced by an earlier sublevel-three node on the branch (as a least-block-element $i$ whose block contains the preimage for some $M$-element $m <_{\mathbb{N}} n$), we don't need the sublevel-one binary split for $n$. Similarly, branching at sublevels two and three for $n$ occurs only if the $M$-element $n$ does not already carry the label of some earlier node on the branch. The {\it on} path through this tree at stage $s$ is defined in the obvious way: the left or right branch of the sublevel-one binary split for each $n \leq_{\mathbb{N}} s$ depending on whether or not $n$ is {\it on} at stage $s$, followed by the leftmost sublevel-two pair $p,q$ that are both {\it on} (if any), followed by the leftmost sublevel-three $i$ that is {\it on} (if any). Note the `if any': the {\it on} branch at stage $s$ will not extend beyond sublevel-one for an $n$ if no $p,q$ for that $M$-element $n$ is {\it on}, and not extend beyond sublevel-two if no suitable $i$ within that $[p,q]$ is {\it on}. Note also that, since the fallow $M$-cluster around $n$ is determined by (the labels from nodes on branches to the left of) the node $p,q$ at which we're acting, we can give a precise definition of the {\it on} path at each stage only later, as part of the description of our construction. For the same reason, the existence of the {\it true} path, the leftmost infinite path whose initial segments will coincide with the {\it on} path infinitely often, can be established only after the description of our construction.

\medskip

So we imagine our $\Pi_{2}$ guessing tree as a subtree of the infinite tree branching downward from the root with three sublevels for each $n = 1, 2, 3, \cdots$: a binary split at sublevel-one for $n$ (from each of the sublevel-three nodes for $n-1$), followed by an $\omega$ split from each of these nodes at sublevel-two, followed by another $\omega$ split from each of these nodes at sublevel-three. Our $\Pi_{2}$ guessing tree will employ a finite, left-most collection of these nodes at each level. If branching at sublevel-one (or sublevels two and three) for an $n$ is unnecessary, we imagine our tree as taking the leftmost branch, with that sublevel-one node referencing no $L$-element $n$ (or the sublevel-two node referencing no $p,q$ and the sublevel-three node no $i$), merely providing a path through to the next level. Otherwise, our guessing tree will include the nodes at sublevels one, two, and three as described in the previous paragraph. We will drive one of the constructions described in the previous paragraphs along each path of this dynamic $\Pi_{2}$ guessing tree, acting along the {\it on} path at each stage. The infinite {\it true} path will be the leftmost one whose initial segments coincide with the {\it on} path infinitely often. Construction along paths to the left of the {\it true} path will eventually cease; construction along paths to the right will always subsequently be rescinded; only the construction along the {\it true} path will drive inexorably forward. The only nodes $\sigma$ at level $n$ that will be used as labels are sublevel-one nodes (to label the image in $M$ of the $L$-block around the least-block-element $n$ referenced by $\sigma$) and sublevel-three nodes (to label the image in $M$ of the $L$-block around the least-block-element $i$ that $\sigma$ references, which image will include the fallow $M$-cluster around $n$). At the end of the construction at each stage the $L$-block around the least-block-element referenced by a node $\sigma$ on the {\it on} path will be in one-to-one correspondence with the $M$-cluster labeled $\sigma$. 

\medskip

The union of these bijective correspondences over all the nodes $\sigma$ along the {\it on} path at stage $s$, between the $L$-block around the element that $\sigma$ references and the $M$-cluster labeled with that $\sigma$, will be used to define the partial function $f^{s}: L^{s} \rightarrow M^{s}$. We have to be a little careful here since, as mentioned before, $L$-elements may belong to more than one $L$-block (and so have more than one image in $M$). $L$-elements referenced by nodes along the {\it on} path will correspond to just one $M$-element at stage $s$. $L$-elements $a$ that are not {\it on} at this stage may belong to as many as two of the $L$-blocks around elements referenced by nodes on the {\it on} path and, consequently, correspond to more than one $M$-element at that stage. If $a$ has two such images, we define $f^{s}(a)$ to be the one defined by the earlier node $\sigma$ (closer to the root). This will make each $f^{s}: L^{s} \rightarrow M^{s}$ a partial bijective function, which functions will allow us to define our $f: L \rightarrow M$ isomorphism: we will establish the existence of a subsequence $s_{i}, s_{2}, s_{3}, \cdots$ of stages on which the construction will resemble a `finite injury' construction of the isomorphism $f$ in that the $f^{s_{i}}$ will eventually include each $L$-element and each $M$-element, and will remain fixed on each of these elements for all but finitely many $i$. As mentioned before, $M$-elements will be ordered when introduced, which order will be preserved, thus making $M$ a computable linear order. Pairs of $M^{s}$-elements not tagged with a common label will never assume a common label at a later stage and will end up in separate blocks in $M$, which ensures that the enumeration at each stage of all pairs of elements in $M^{s}$ that do not have a common label will provide a computable enumeration of $M$'s non-block relation $\neg B_{M}(x,y)$. 

\medskip

A delicacy in our setup, crucial to its success, bears highlighting: The fallow $M$-cluster around $n$ and its selected preimage within an interval $[p,q]$ are determined at the sublevel-two node referencing that $p,q$ (in step 5 of the construction), not at the sublevel-three nodes below, at which we will act (in step 8) to incorporate $n$ and its fallow $M$-cluster into the image of an $L$-block around an {\it on} element $i$ within $[p,q]$. That fallow $M$-cluster around $n$ is defined by the labels that $n$ carries from nodes on paths to the left of this node $p,q$ ... whereas, when acting at the sublevel-three node $i$ below this node $p,q$, we must in fact incorporate into the image of the $L$-block around $i$ the larger fallow $M$-cluster around $n$ defined by {\it all} labels that $n$ carries from nodes on paths to the left of the node $i$, which is one sublevel below the node $p,q$. The additional labels on $n$ that make this fallow $M$-cluster larger would come (only) from the sublevel-three nodes $j$ immediately below this node $p,q$ and to the left of the node $i$ at which we're acting, labels which were placed to mark the images of the $L$-blocks around those $j$. We are able to skirt this obstacle by dint of the fact that, if $i$ is indeed the true least-block-element whose block contains the selected preimage of the fallow $M$-cluster around $n$ as calculated at the node $p,q$, then, eventually, those $j$ {\it and their $L$-blocks} will all be contained within the $L$-block around $i$.  

\bigskip

\noindent
{\bf {\it The Construction}}: 

\medskip

For each $s = 1, 2, 3, \cdots$ in turn begin the construction at stage $s$ by starting at the root of our $\Pi_{2}$ guessing tree and working downward, defining the {\it on} path through the three sublevels for each $n \leq_{\mathbb{N}}s$, and acting as described below at each node along that path. Our action will define for some of those nodes an $L$-element that the node references and a cluster of $M$-elements labeled with that node, which labeled cluster will be in bijective correspondence with the block around the referenced $L$-element. The order $<_{M}$ is inherited from $<_{L}$ and defined as elements are introduced into $M$. At the end of the construction at each stage, all pairs of $M$-elements that do not have a common label are enumerated into $M$'s non-block relation $\neg B_{M}(x,y)$, and the bijective correspondences between the $L$-elements referenced by nodes along the {\it on} path and the $M$-clusters labeled with those nodes are used as described in the preliminary definitions section to define the partial bijection $f^{s}: L^{s} \rightarrow M^{s}$ (each $f^{s}(a)$ is defined to be the image in the bijective correspondence defined by the earliest node on the {\it on} path whose $L$-block includes $a$). One piece of the algorithm is executed repeatedly through the construction at each stage and is presented here as a blanket instruction: As each node along the {\it on} path is defined, delete the effects of all earlier action at nodes on all paths to the right (i.e.\ consider those nodes as no longer referencing any $L$-elements and remove all labels in $M$ tagged with those nodes). 

\medskip

Having completed the construction at stage $s$ for $1, 2, 3, \cdots, n-1$, continue the construction for $n$ ($\leq_{\mathbb{N}} s$) by performing in order the following steps: 

\bigskip

At sublevel-one for $n$: 

\begin{itemize}

\item[1:] If $n$ is not {\it on} at stage $s$ the {\it on} path takes the right branch at sublevel-one and no further action is performed at this sublevel; go to step 3. Otherwise the {\it on} path takes the left branch. If $n$ is {\it on} but is already referenced by an earlier node along the {\it on} path, then no further action is performed at this sublevel; go to step 3. Otherwise $n$ is referenced by this node at sublevel-one for $n$.

\item[2:] Arrange that the $M$-cluster labeled with this node $n$ is in bijective correspondence with the $L$-block around $n$, with the correspondence mapping $n$ to the $\mathbb{N}$-least element in the $M$-cluster. Begin with the existing $M$-cluster around $n$ labeled with this node; build the cluster from scratch if none exists. Remove labels from some of the outermost elements in the cluster, or add new elements labeled with this node to the outer ends of the cluster (using the $\mathbb{N}$-least elements that have not yet been enumerated into $M$), as necessary, to make the $M$-cluster the same size as the $L$-block around $n$ at this stage. Do this carefully: remove labels from elements or add new elements to the left end of the $M$-cluster so that there are exactly as many elements in the cluster to the left of its $\mathbb{N}$-least element as there are to the left of $n$ in the $L$-block around $n$; then handle the right side similarly. If this is a completely new $M$-cluster situate it correctly in $M$, with order inherited from $L$, with respect to all $M$-clusters that are labeled with earlier nodes along the {\it on} path (and as far left as possible with respect to all other $M$-elements).

\end{itemize}

At sublevel-two for $n$: 

\begin{itemize} 

\item[3:] If the $M$-element $n$ is labeled with an earlier node along the current {\it on} path, then the {\it on} path goes through the leftmost nodes at sublevels two and three for $n$ and no further action is performed at these sublevels. Go to step 1 of the construction for $n+1$. 

\item[4:] Otherwise, consider the location in $M$ of this element $n$ with respect to all $M$-clusters labeled with earlier nodes along the current {\it on} path. We must find a preimage for $n$ that lies in the same interval in $L$ with respect to the preimages of these higher priority labeled $M$-clusters. The sublevel-two nodes for $n$ are taken to reference, from left to right, a lexicographic ordering of all pairs $p,q$ that lie within this interval in $L^{s}$. If none of these pairs $p,q$ is {\it on} (i.e.\ both $p$ and $q$ {\it on} at stage $s$), the {\it on} path does not extend to sublevel-two for $n$ and no further action is performed at these sublevels; go to step 1 of the construction for $n+1$. If one of these pairs $p,q$ is {\it on}, then the next node along the {\it on} path is the leftmost of these $p,q$ that is {\it on} at this stage. 

\item[5:] Consider the selected preimage within $[p,q]$ (the leftmost in the lexicographic list of contiguous sequences of the right size) for the fallow $M$-cluster around $n$ (defined by the labels on $n$ from all nodes on paths to the left of $p,q$), both of which are described in the preliminary definitions section. If too few elements have been enumerated into $[p,q]$ for a preimage to exist, then the {\it on} path does not extend beyond sublevel-two for $n$ and no further action is performed at these sublevels; go to step 1 of the construction for $n+1$. If this fallow $M$-cluster around $n$ or its selected preimage within this $[p,q]$ has changed since the last time the node $p,q$ was on the {\it on} path, delete the effect of all earlier action at nodes on all paths below this node. 

\end{itemize}

At sublevel-three for $n$: 

\begin{itemize} 

\item[6:] Consider the leftmost sublevel-three nodes immediately below the node $p,q$ as referencing a list of all $i$ in $[p,q]$ that are $\leq_{\mathbb{N}}$ the ${\mathbb{N}}$-least element of the selected preimage (of step 5), ordered from left to right according to the last time that an element was enumerated into $L$ between $i$ and the preimage (as described in the preliminary definitions section). If the $i$ referenced by one of these nodes has changed since the last stage, delete the effect of all earlier action at that node and at all nodes on paths below it. For the remaining sublevel-three nodes (that do reference an $i$), cut back, if necessary, the $M$-cluster labeled with that node so that it is no larger than the $L$-block around the $i$ that the node references (i.e.\ remove that node's label from all $M$-elements in the cluster that correspond to $L$-elements that were previously in the $L$-block around $i$ but are no longer). 

\item[7:] If there is no $i$ on this list that is {\it on} at this stage and whose current block in $L$ wholly contains the selected preimage of the fallow $M$-cluster around $n$ {\it and} the block around the $L$-element referenced by each sublevel-three node below this $p,q$ and to the left of the node $i$, then the {\it on} path does not extend to sublevel-three for $n$ and no further action is performed at this sublevel; go to step 1 of the construction for $n+1$. Otherwise the next node along the {\it on} path is the leftmost node $i$ with these properties.

\item[8:] Arrange that the $M$-cluster labeled with this node $i$ (includes the {\it whole} labeled block around $n$ and) is in bijective correspondence with the $L$-block around $i$, which correspondence maps the selected preimage to the fallow $M$-cluster around $n$. Begin with the $M$-cluster around $n$ that includes not only the fallow $M$-cluster around $n$ as calculated at the node $p,q$ (see step 5) but also all other labels on $n$ (from all sublevel-three nodes below $p,q$ and on or to the left of the {\it on} path). Add, if necessary, new elements to the outer ends (using the $\mathbb{N}$-least elements that have not yet been enumerated into $M$) to make this $M$-cluster around $n$ (labeled with the node $i$) the same size as the $L$-block around $i$. Note that $i$ was selected (in step 7) to make this possible. Do this carefully: begin with the $M$-cluster around $n$ consisting of all elements that share a label with $n$ and add, if necessary, new elements to the left end of the cluster so that there are exactly as many elements in the cluster to the left of the fallow $M$-cluster around $n$ as there are to the left of its selected preimage in the $L$-block around $i$; then handle the right side similarly. This ends the construction for $n$; go to step 1 of the construction for $n+1$.

\end{itemize} 

\medskip

\noindent
{\bf {\it Final Verification}}: 

\medskip

Preliminary observations: 
\begin{itemize}
\item $M$ is a computable linear order: Elements are ordered as they are introduced into $M$, which order is never changed. 
\item $M$'s non-block relation $\neg B_{M}(x,y)$ is computably enumerable: Elements in $M$ that have a common label will never increase their distance from each other so long as they retain that label (since no elements are ever introduced internally into a labeled block), and elements that do not have a common label will never assume one at a later stage (only new elements are used to grow $M$-clusters). 
\item All $M$-clusters with labels from nodes along the {\it on} path will be ordered with respect to each other exactly as their preimages are ordered in $L$, and no $M$-element will ever belong to more than one of these $M$-clusters (labeled with nodes along the {\it on} path). 
\item A labeled $M$-cluster will grow only at stages when the labeling node is on the {\it on} path and will, at the end of those stages, be the same size as the block at that stage around the $L$-element that the node references. 
\item An $M$-cluster labeled with a node will disappear completely (all labels tagged with that node removed) only if the node is on a path to the right of the {\it on} path, or below a sublevel-two node $p,q$ at which the fallow $M$-cluster around $n$ or its selected preimage have been redefined (step 5), or at or below a sublevel-three node whose referent has changed (step 6). 
\end{itemize}

\medskip

We show now that there is an infinite {\it true} path which, for every level $n$, is the leftmost path down to that level that coincides with the {\it on} path infinitely often; that every node on this path will change the $L$-element that it references only finitely often, with its final referents (if any) being true least-block-elements; and that every $M$-cluster labeled with a node on this path will disappear (completely) only finitely many times. We will establish these facts simultaneously, by induction on $n$, down to the three sublevels for each $n$. 

\medskip

Assume that we have established the abovementioned facts for $1, 2, \cdots n-1$ and let $s$ be a stage after which the {\it on} branch is never to the left of the {\it true} branch all the way down to sublevel-three for $n-1$; by which stage all these nodes have taken on their final referents; and after which stage no $M$-cluster labeled with one of these nodes ever (completely) disappears. We will now establish the facts for the three sublevels for $n$. All stages mentioned below are intended to be subsequent to $s$.  
\begin{itemize}
\item At sublevel-one for $n$: If $n$ is not a true least-block-element consider a stage after which $n$ never again appears {\it on}. At the next stage when the {\it on} path coincides with the {\it true} path down to sublevel-three for $n-1$, the {\it on} path will take the right branch at sublevel-one for $n$ as will the {\it true} path: the {\it on} path will never again be left of this node; this node will never again take on a referent; nor will it ever again label an $M$-cluster. If $n$ is a true least-block-element, consider the next stage when $n$ is {\it on} and the {\it on} path coincides with the {\it true} path down to sublevel-three for $n-1$; the {\it on} path will at this stage take the left branch at sublevel-one for $n$ as will the {\it true} path, and the {\it on} path will never again be left of this node. If $n$ is already referenced by some (sublevel-three) node earlier on the {\it true} path, this node will never again take on a referent nor will it ever again label an $M$-cluster. If $n$ is not referenced by any earlier node on the {\it true} path, then it will subsequently always be the referent of this node, which will label an $M$-cluster (corresponding to the $L$-block around $n$) that will never (completely) disappear.
\item At sublevel-two for $n$: Consider a subsequent stage by which the $M$-element $n$ has assumed its final position (in or out) with respect to all $M$-clusters labeled with nodes down to sublevel-one for $n$ on the {\it true} path ($M$-elements that do not share a label will never subsequently assume a common one), and when the {\it on} path coincides with the {\it true} path down to sublevel-one for $n$. If $n$ belongs to (precisely) one of the $M$-clusters labeled with an earlier node on the {\it true} path then, at this stage, the {\it on} path (and the {\it true} path) will take the the leftmost branch at sublevels two and three for $n$; the {\it on} path will never again be left of these nodes;  neither of these nodes will ever again take on a referent; nor will they ever again label an $M$-cluster. If $n$ belongs to no $M$-cluster labeled with an earlier node on the {\it true} path, consider the interval in $L$ between the elements referenced by these earlier nodes on the {\it true} path between whose images $n$ lies. Since this interval is infinite (the $L$-elements referenced by these nodes are all true least-block-elements), there will be a lexicographically least pair $p,q$ of elements within this $L$-interval both of which are true least-block-elements. Consider a subsequent stage when $p,q$ is {\it on} (i.e.\ both $p$ and $q$ are {\it on}) and after which no lexicographically earlier pair within this $L$-interval ever appears {\it on}. At this stage the fallow $M$-cluster around $n$ as calculated at $p,q$ will have taken on its final form (no further labels from nodes on paths to the left of $p,q$ will be added). After this stage, every time the {\it on} path coincides with the {\it true} path down to sublevel-one for $n$, it will take the branch through this node $p,q$ at sublevel-two for $n$. Consider a further stage when the selected preimage within this interval $[p,q]$ for the fallow $M$-block around $n$ (the leftmost on the list of candidates) is indeed a true set of contiguous $L$-elements. 
\item At sublevel-three for $n$: Let $i$ be the true least-block-element within this $L$-interval $[p,q]$ whose (true) block contains the abovementioned selected preimage for the fallow $M$-block around $n$. Consider a subsequent stage when the (finite) contiguous part of the $L$-block around the selected preimage that contains all elements that are $\leq_{\mathbb{N}}$ than the $\mathbb{N}$-least element of the selected preimage have been enumerated into $L$ and after which none of these elements (other than $i$) are ever again {\it on}. After this stage the leftmost sublevel-three nodes below node $p,q$ will refer to these elements (that are $\mathbb{N}$-less than the $\mathbb{N}$-least element of the selected preimage), one of which will be $i$, and there will be no redefining of referents at these nodes. Consider a subsequent stage when $i$ is {\it on} and when the {\it on} path coincides with the {\it true} path down to sublevel-two for $n$; the {\it on} path will at this stage take the branch through node $i$ at sublevel-three for $n$ as will the {\it true} path; the {\it on} path will never again be left of this node; this node will always refer to $i$; and there will be an $M$-cluster labeled with this node $i$ that will always contain the fallow $M$-cluster around $n$.
\end{itemize}

This ends the proof of the abovementioned three facts. 

\medskip

Consider the subsequence $s_{1}, s_{2}, s_{3}, \cdots$ of stages in which each $s_{n}$ is the first stage after $s_{n-1}$ at which the {\it on} path coincides with the {\it true} path down to sublevel-three for $n$, after which the {\it on} path is never to the left of these nodes, after which the $L$-elements that these nodes reference are never changed, and after which the $M$-clusters labeled with these nodes are never (completely) deleted. If, at this stage $s_{n}$, the node $\sigma$ on the {\it true} path at sublevel-one for $n$ references $n$, then, since the bijection defined at $\sigma$ between the $L$-block around $n$ and its image (the $M$-cluster labeled with $\sigma$) always maps $n$ to the $\mathbb{N}$-least element of the image (see step 2), we have that $f^{s_{n}}(n) = \lim_{i \rightarrow \infty}f^{s_{i}}(n)$. Similarly, if, at this stage $s_{n}$, the node $\sigma$ on the {\it true} path at sublevel-three for $n$ references an $L$-element $i$ (whose $L$-block contains the selected preimage for the fallow $M$-cluster around $n$), then, since the bijection defined at that $\sigma$ is centered around a map of the selected preimage onto the fallow $M$-cluster (see step 8), it follows that the $f^{s_{i}}$ have by stage $s_{n}$ taken on their final limiting value on the $M$-element $n$. Moreover, since no more elements will be enumerated into $L$ between $i$ and the selected preimage after this stage, these functions have by this stage also taken on their final limiting value on the least-block-element $i$. So the partial bijections $f^{s_{i}}$ eventually take on final limiting values on every $M$-element and on every least-block-element in $L$. 

\medskip

For the other elements of $L$, the non least-block-elements $a$, observe that, once all the elements between $a$ and the least-block-element $n$ of $a$'s block have been enumerated into $L$ and have ceased coming {\it on}, the image of $a$ in the bijection between the $L$-block around $n$ and its image will remain fixed. The only complication is that $a$ may appear to be in the $L$-blocks around more than one of the elements referenced by nodes on the {\it true} path and $f^{s_{i}}(a)$ may be defined via the $L$-block around the element referenced by some node earlier than that which references $n$ (an erroneous definition in that $a$ is not truly a member of that earlier $L$-block). But, once the {\it on} path coincides at these stages $s_{i}$ with the {\it true} path far enough down to include a pair of true least-block-elements around $n$ that separates it from the elements referenced by all earlier nodes on the {\it true} path, then, at all subsequent stages, the $f^{s_{i}}(a)$ will be defined correctly, via the correspondence between the $L$-block around $n$ and its image. 

\medskip

Thus $f(x) = \lim_{i \rightarrow \infty}f^{s_{i}}(x)$ over these stages $s_{i}$ defines an isomorphism from the given computable linear order $L$ to the constructed $M$, which is a computable linear order with computably enumerable non-block relation $\neg B_{M}(x,y)$.

\hspace{\fill} $\Box$

\medskip

Extending the theorem: Our construction requires that the condensation-type of $L$ is not just dense but is dense without endpoints (i.e.\ $\eta$ rather than $1 + \eta$, $\eta + 1$, or $1 + \eta + 1$); the result carries over nonetheless ...

\medskip

\noindent
{\bf THEOREM 1 EXTENDED:}

\medskip

{\it Every computable linear order with dense condensation-type and with no infinite, strongly $\eta$-like interval has a computable copy with computably enumerable non-block relation.}

\medskip

\noindent
{\bf PROOF:}

\medskip

Consider any computable linear order $L$ with condensation-type $1 + \eta$ and with no infinite, strongly $\eta$-like interval. Form the computable linear order $L^{\prime}$ by adding to the left of $L$ a decidable copy of $\omega \times \eta$ to produce a computable linear order $L^{\prime}$ with condensation-type $\eta$. Apply Theorem 1 to $L^{\prime}$ to produce a computable copy $M^{\prime}$ with computably enumerable non-block relation. Produce from this $M^{\prime}$ an $M \cong L$ by deleting everything in $M^{\prime}$ to the left of the element that corresponds to the leftmost point of $L$. If $L$ has no leftmost point, i.e.\ if its leftmost block is of type $\omega^{\star}$ or $\omega^{\star} + \omega$, then delete everything in $M^{\prime}$ to the left of an element that corresponds to some point in this leftmost block of $L$, and add to the left of the resulting linear order a decidable copy of $\omega^{\star}$. In either case, this surgery will produce a computable copy $M$ of $L$ with computably enumerable non-block relation. Computable $L$ with condensation-types $\eta +1$ and $1 + \eta + 1$ are handled similarly.  

\hspace{\fill} $\Box$

\bigskip

Exceptions that probe the rule: Crucial use is made in our construction (in step 4) of the existence of the pair of least-block-elements $p,q$ in $L$ between least-block-elements of higher priority. Their existence could not be guaranteed if $L$'s condensation-type were not $\eta$; nor does the result carry over to such linear orders: {\it There is a computable linear order with no infinite, strongly $\eta$-like interval that has no computable copy with computably enumerable non-block relation}. (A {\it discrete} example of such a linear order, i.e.\ with every block of type $\omega^{\star} + \omega$, is presented in Moses [1988].) 

\medskip

Similarly, crucial use is made (in step 5) of the existence of arbitrarily large blocks within every infinite $L$-interval, which could not be guaranteed if $L$ had an infinite strongly $\eta$-like interval; nor does the result carry over to such linear orders: {\it There is a computable linear order with dense condensation-type that has no computable copy with computably enumerable non-block relation.} (In [1997] Coles, Downey, and Khoussainov present an example of such a linear order, with condensation-type $\eta + 1$ which in fact has no computable copy in which the initial segment with condensation-type $\eta$, i.e.\ everything but the final $\omega^{\star}$-block, is computably enumerable.) 

\medskip

And, every time its node comes {\it on}, a labeled $M$-block shakes off `outer' elements, labeled with nodes on branches to the right, enumerating some pairs of $M$-elements that were previously considered to be in the same block into the non-block relation $\neg B_{M}(x,y)$, which means that we cannot guarantee that the complement $B_{M}(x,y)$ is computably enumerable. Nor can the result be strengthened to guarantee this: {\it There is a computable linear order with dense condensation-type and no infinite, strongly $\eta$-like interval no computable copy of which has its block relation $B(x,y)$ computably enumerable.} To construct such a linear order begin with a $\Pi_{2}$ suborder of the standard computable linear order of type $\eta$ that has no $\Pi_{1}$ copy. (A $\Pi_{n}$ linear order is that defined as a suborder of the standard computable linear of type $\eta$, or indeed of any computable linear order, by a $\Pi_{n}$ subset of the universe $\mathbb{N}$. That there is, for each $n$, a $\Pi_{n+1}$ linear order that has no $\Pi_{n}$ copy is well known, see Rosenstein [1982].) Produce a computable linear order $L$ from this computable copy of $\eta$ by replacing each element $i$ with a block of size $p_{i}$, the $i$th prime, if $i$ is an element of the $\Pi_{2}$ suborder, and a block of type $\omega^{\star}+\omega$ if it is not. (That such a computable linear order can be produced from any $\Pi_{2}$ linear order is established in Fellner [1976].) $L$ has condensation-type $\eta$ and, since no two finite blocks are of the same size, has no infinite, strongly $\eta$-like interval. And it has no computable copy $M$ with computable block relation since, if it did, the set $\{ x : \forall y >_{M} x \: \neg B_{M}(x,y)\}$ defines a $\Pi_{1}$ copy of the $\Pi_{2}$ linear order that we began with, which does not exist.

\newpage
\begin{center}
{\bf Non-Trivial Self-Embeddings}
\end{center}

\bigskip

In [1940] Dushnik and Miller observed that every countable linear order has a non-trivial self-embedding. It is well known that this does not carry over to computable linear orders: not every computable linear order has a computable non-trivial self-embedding. A long-standing conjecture (see [1987] Kierstead) seeks to characterize the computable linear orders every computable copy of which has a computable non-trivial self-embedding as precisely those with an infinite, strongly $\eta$-like interval. The conjecture remains open. Our Theorem 1 allows us to establish the conjecture for all linear orders with dense condensation-type (see Corollary 2 below). This is a slightly stronger result than the main theorem of Downey, Kastermans, and Lempp [2009], which establishes this result for all $\eta$-like linear orders (which have dense condensation-type but no infinite blocks). 

\medskip

Surprisingly, at least to me, our theorem allows us to characterize completely the other side of the coin: 

\bigskip

\noindent
{\bf COROLLARY 1:}

\medskip

{\it Every computable linear order has a computable copy with a computable non-trivial self-embedding.} 

\medskip

\noindent
{\bf PROOF:}

\medskip

If the linear order has a pair of adjacent blocks then it has an interval of type $\omega + 1$ or $1 + \omega^{\star}$ or $\omega^{\star} + \omega$. Create a computable copy of the linear order by replacing this closed interval with a decidable copy of the same type (and leaving everything else the same). This computable copy has a computable non-trivial self-embedding (the identity everywhere outside the interval and, inside it, mapping each element to its immediate successor in the $\omega$-block and its immediate predecessor in the $\omega^{\star}$-block, both of which are computably identifiable). If the linear order has an infinite, strongly $\eta$-like interval then it clearly has a computable non-trivial self-embedding (mapping elements within that interval into separate blocks). By our Extended Theorem 1, the remaining computable linear orders all have computable copies with computably enumerable non-block relation, which, as in the strongly $\eta$-like interval case, can be used to define a computable non-trivial self-embedding of that copy by mapping elements into separate blocks.

\hspace{\fill} $\Box$

\bigskip

\noindent
{\bf COROLLARY 2:} 

\medskip

The self-embedding conjecture holds for linear orders with dense condensation-type: {\it Every computable copy of a computable linear order with dense condensation-type has a computable non-trivial self-embedding if and only if the linear order contains an infinite, strongly $\eta$-like interval.}  

\medskip

\noindent
{\bf PROOF SKETCH:} 

\medskip

Consider the straightforward technique for constructing, from a computable linear order $L$, a computable copy $M$ with no computable non-trivial self-embedding: a `finite-injury' construction, meshing the enumeration of a computable copy with the diagonalization across all partial computable functions, the candidates for a possible non-trivial self-embedding. The requirements for $n$, with priority decreasing as $n$ increases, would arrange for an image in $M$ for $n \in L$, for a preimage in $L$ for $n \in M$, and that the $n$th partial computable function $f_{n}$ is not a self-embedding on $M$, with, as usual, action to meet the requirements for $n$ being conducted without redefining the $L \rightarrow M$ isomorphism on the points that witness our meeting the requirements with priority higher than $n$. An obvious strategy for ensuring that $f_{n}$ is not a self-embedding on $M$, used for instance to similar purpose in Moses [1988] (on discrete linear orders), is to wait until images $f_{n}(a)$ and $f_{n}^{2}(a)$ are defined for some $a \in M$ with all three within the same interval between points of higher priority (on which the $L \rightarrow M$ isomorphism cannot be redefined) and with either $a <_{M} f_{n}(a) <_{M} f_{n}^{2}(a)$ or $a >_{M} f_{n}(a) >_{M} f_{n}^{2}(a)$. Then, every time a new point is enumerated into this interval (between higher priority points) in $L$, redefine the $L \rightarrow M$ isomorphism within this higher-priority interval to `feed' the $M$-interval $[ a, f_{n}(a)]$ until it is larger than $[ f_{n}(a), f_{n}^{2}(a)]$. Once this has been achieved, all that remains is to ensure that the $M$-interval $[ f_{n}(a), f_{n}^{2}(a)]$ stops growing, which can be arranged by identifying a contiguous set of $L$-elements of the right size within this higher-priority interval and making it the preimage of the $M$-interval $[ f_{n}(a), f_{n}^{2}(a)]$. Since every infinite $L$-interval contains arbitrarily large blocks, such a contiguous set of $L$-elements must exist, and we can identify one by maintaining a lexicographic list of all contiguous sets of the right size, deleting those that cease being contiguous (as we did in selecting a preimage for the fallow $M$-block in the proof of Theorem 1), which will eventually leave a true contiguous set as the leftmost in the list. The problem is that we may not be able to redefine the $L \rightarrow M$ sufficiently so as to make this selected contiguous sequence the preimage of $[ f_{n}(a), f_{n}^{2}(a)]$, which would (could) happen if the selected preimage is in fact in the same block as one of the higher-priority elements. We got around this in our proof of Theorem 1 by using the points $p,q$, which ensured that the selected preimage was in a separate block from those of the higher-priority elements. (In Moses [1988], which dealt with discrete linear orders, we just kept moving the preimage of $[ f_{n}(a), f_{n}^{2}(a)]$ toward one of the higher-priority elements, which eventually left it in the (infinite) block around that element.) If $L$ has its non-block relation computably enumerable, we can use it to select a preimage for $[ f_{n}(a), f_{n}^{2}(a)]$ that is in a separate block from all higher-priority elements, and so carry through the construction. What this establishes is stronger than Corollary 2: {\it The self-embedding conjecture holds for all linear orders that have a computable copy with computably enumerable non-block relation.}

\hspace{\fill} $\Box$

\newpage
\begin{center}
{\bf References}
\end{center}

\begin{itemize}

\item[][1981] C.\ Ash and A.\ Nerode, {\it Intrinsically recursive relations}, `Aspects of effective algebra', edited by J.N.\ Crossley (Upside Down A Book Company, Steel's Creek, Australia), 26-41

\item[][1997] R.\ Coles, R.\ Downey, and B.\ Khoussainov, {\it On initial segments of computable linear orders}, Order {\bf 14}, 107-124

\item[][2009] R.\ Downey, B.\ Kastermans, and S.\ Lempp, {\it On computable self-embeddings of computable linear orders}, to appear

\item[][1940] B.\ Dushnik and and E.\ Miller, {\it Concerning similarity transformations of linearly ordered sets}, Bull.\ Amer.\ Math.\ Soc.\ {\bf 46}, 322-326

\item[][1970] L.\ Feiner, {\it Hierarchies of Boolean algebras}, J.\ Symbolic Logic {\bf 35}, 365-373 

\item[][1976] S.\ Fellner, {\it Recursiveness and finite axiomatizability of linear orderings}, Thesis Rutgers University, N.J.

\item[][1966] C.\ Jockusch, Jr., {\it Reducibilities in recursive function theory}, Thesis MIT, Cambridge, Mass.

\item[][1991] C.\ Jockusch, Jr.\ and R.\ Soare, {\it Degrees of orderings not isomorphic to recursive linear orderings}, `International Symposium on Mathematical Logic and its Applications (Nagoya, 1988)', Ann.\ Pure Appl.\ Logic {\bf 52}, 39-64 

\item[][1987] H.\ Kierstead, {\it On $\Pi_{1}$ automorphisms of recursive linear orders}, J.\ Symbolic Logic {\bf 52}, 681-688 

\item[][1986] M.\ Moses, {\it Relations intrinsically recursive in linear orders}, Zeitschr.\ f.\ math.\ Logik und Grundlagen d.\ Math.\ {\bf 32}, 467-472 

\item[][1988] M.\ Moses, {Decidable discrete linear orders}, J.\ Symbolic Logic {\bf 53}, 531-539

\item[][1982] J.\ Rosenstein, {\it Linear orderings}, Academic Press, New York/London 

\end{itemize}

\end{large}

\end{document}